\documentclass[12pt, a4paper]{article}
\usepackage{a4wide}
\usepackage{amsmath}
\usepackage{amsfonts}
\usepackage{amssymb}
\usepackage{graphicx}
\usepackage{xcolor}
\graphicspath{{./Figures/}}

\DeclareMathOperator{\lip}{lip}
\DeclareMathOperator{\dom}{dom}

\DeclareMathOperator{\gph}{gph}
\DeclareMathOperator{\clm}{clm}
\DeclareMathOperator{\Lipusc}{Lipusc}

\DeclareMathOperator{\cone}{cone}

\DeclareMathOperator{\Limsup}{Lim sup}

\newtheorem{theo}{Theorem}

\newtheorem{prop}{Proposition}
\newtheorem{cor}{Corollary}
\newtheorem{rem}{Remark}
\newtheorem{exa}{Example}

\newenvironment{dem}[1][Proof]{\noindent \textbf{#1.} }{\ \rule{0.5em}{0.5em}}

\begin{document}

\title{Bridging local and semilocal stability: A topological approach
\thanks{
This research was partially supported by the Spanish Ministry of Science and Innovation (MICINN) through grant PID2022-136399NB-C22, and by the European Regional Development Fund (ERDF), ``A way to make Europe.''
}}
\author{J. Camacho
\thanks{
Center of Operations Research, Miguel Hern\'{a}ndez University of Elche, 03202 Elche (Alicante), Spain (j.camacho@umh.es).}}
\date{}
\maketitle

\begin{abstract}
This paper establishes a general topological condition under which the semilocal stability of a set-valued mapping can be exactly determined by its local stability properties. Specifically, we investigate the relationship between the Lipschitz upper semicontinuity modulus ---a semilocal measure of variation for the image set--- and the local calmness moduli. While these two quantities coincide for mappings with convex graphs, the relationship generally breaks down in the absence of convexity, making the semilocal modulus exceptionally difficult to compute. We prove that if a mapping is outer semicontinuous in the Painlev\'{e}-Kuratowski sense and locally compact around the nominal parameter, the Lipschitz upper semicontinuity modulus is exactly the supremum of the local calmness moduli over the nominal set. In addition to the theoretical advance, this equality enables the precise calculation of semilocal error bounds via point-based formulae. We illustrate the broad applicability of this theorem by setting it up in several non-convex frameworks in parametric optimization, including piecewise convex and semi-algebraic mappings, feasible and optimal set mappings under full data perturbations, generalized equations and linear complementarity problems, semi-infinite inequality systems, and parameterized sub-level sets.


\end{abstract}

\section{Introduction}\label{sec: intro}    

Variational analysis is a cornerstone of mathematical optimization, providing essential insights into the behavior of the image of a set-valued mapping under data perturbations. A central challenge in this field is quantifying the rate at which such images (usually refereed to as solutions) expand relative to perturbations in its parameters. The quest to quantify this stability in the finite linear inequality system setup dates back to the seminal work of Hoffman \cite{hoffman1952approximate}, who established the existence of a global error bound. This paper focuses on the property of \emph{Lipschitz upper semicontinuity} and its associated modulus, a semilocal quantifier that captures the behavior of the entire image set of a generalized mapping near a nominal parameter. 

It is important to note that the terminology regarding this property is not entirely standard in the literature. While the underlying concept relates closely to the \emph{upper Lipschitz continuity} discussed in the classic works of Robinson \cite{robinson1981some} and the \emph{outer Lipschitz continuity} properties found in Rockafellar and Wets \cite{rockafellar1998variational} or Dontchev and Rockafellar \cite{dontchev2009implicit}, we adopt the term \emph{Lipschitz upper semicontinuity} following the nomenclature used by Klatte and Kummer \cite{klatte2002nonsmooth} or Uderzo \cite{uderzo2021quantitative}. 

Stability results of this type have profound applications across numerous domains, acting as the theoretical backbone for algorithms dealing with uncertainty and hierarchical decision-making. For instance, the classical works of Fiacco \cite{fiacco1983introduction} and Bonnans and Shapiro \cite{bonnans2000perturbation} on sensitivity analysis in nonlinear programming extensively rely on such bounds to ensure well-behaved solution trajectories under data variations. In the realm of stochastic programming, Shapiro, Dentcheva, and Ruszczyński \cite{shapiro2009lectures} demonstrate that qualitative and quantitative stability under probability measure perturbations is strictly governed by Lipschitz-type properties of the corresponding optimal set mappings. Furthermore, these stability conditions are indispensable in analyzing hierarchical models; as established by Luo, Pang, and Ralph \cite{luo1996mathematical} for mathematical programs with equilibrium constraints (MPECs) and by Dempe \cite{dempe2002foundations} for bilevel optimization problems, they are fundamental to establishing convergence and deriving optimality conditions.

The analysis of stability is typically conducted from two distinct perspectives: local and semilocal (or global). The \emph{calmness} property is a local notion, measuring the behavior of the mapping at a specific point in the graph restricted to a neighborhood. In contrast, the Lipschitz upper semicontinuity property considers the semilocal behavior, taking into account variations over the entire set of nominal images. A recurring fundamental question is whether the semilocal stability is determined solely by local information. That is, does the Lipschitz upper semicontinuity property holds whenever the calmness property does across the nominal image? If so, does the Lipschitz upper semicontinuity modulus strictly coincide with the supremum of the calmness moduli? 

While calmness moduli are point-based and often explicitly computable via Karush-Kuhn-Tucker (KKT) systems, coderivatives, or generalized differentiation tools --as developed in the foundational text of Mordukhovich \cite{mordukhovich2006variational}, advanced in the linear context within multiple collaborations by Cánovas, Parra and others \cite{ canovas2014calmness3,canovas2016calmness, canovas2016outer,canovas2014calmness2,canovas2014calmness, gisbert2018calmness}, and extended to the non-convex settings by Gfrerer \cite{gfrerer2013first} or Henrion and Jourani \cite{henrion2001subdifferential}-- the Lipschitz upper semicontinuity modulus is a set-based supremum that is generally much harder to estimate. For mappings with convex graphs, a classic result by Robinson \cite{robinson1981some} establishes that the global Lipschitz constant exactly coincides with the supremum of the local calmness moduli. This equality is a powerful theoretical tool, as it allows one to compute a global stability measure simply by investigating local properties. However, outside the convex setting, a priori, this relationship might break down, and the Lipschitz upper semicontinuity modulus can be strictly greater than the supremum of calmness moduli. 

In this paper, we address the challenge of identifying a broader class of mappings for which the Lipschitz upper semicontinuity modulus is exactly the supremum of the calmness moduli. While previous works have focused on specific structures (see, e.g., \cite{camacho2022calmness, camacho2023lipschitz, camacho2026lipschitz}) such as convex-graph mappings or polyhedral ones, we present a general topological sufficient condition. Concerning the Lipschitz upper semicontinuity property itself, the immediate antecedent can be traced to Klatte \cite{klatte2026lipschitz}, where the finiteness of the corresponding modulus is characterized in the context of quadratic problems restricted to polyhedral feasible sets. 

Our main result, Theorem \ref{theo closed outer}, establishes that outer semicontinuity of the mapping, combined with local compactness around the nominal parameter, is sufficient to guarantee this equality. This result significantly extends the scope of applicability, covering many non-convex mappings encountered in practice. Indeed, as a straightforward consequence, any locally compact mapping with a closed graph yields the desired equality, as Corollary \ref{cor closed outer} shows.

Furthermore, we explore the implications of this result for several well-known classes of mappings in optimization theory: namely, the feasible and argmin mappings under full perturbations, linear complementarity problems, general polyhedral feasible sets, and semi-infinite inequality systems. These structures, which naturally arise in parametric linear programs, lack convexity in their graphs but satisfy the proposed topological conditions, thereby allowing their semilocal stability to be fully characterized by local calmness moduli.

The paper is organized as follows. Section \ref{sec: preliminaries} introduces the necessary notation and preliminaries on variational analysis, including formal definitions of the moduli. Section \ref{sec: main} presents the main theoretical result, proving the equality of moduli for outer semicontinuous mappings with compact nominal sets. Section \ref{sec: classes of mappings} recalls the class of wcpc mappings and applies the main theorem to this context, together with the optimal and feasible set mappings in the full perturbation framework, linear complementarity problems, general polyhedral feasible sets, and semi-infinite inequality systems. Finally, Section \ref{sec: conclusions} provides concluding remarks and perspectives for future research.

\section{Preliminaries}\label{sec: preliminaries}

Throughout this paper, let $X$ and $Y$ be metric spaces equipped with distances $d_X$ and $d_Y$, respectively. For a point $x \in X$ and a subset $\Omega \subset X$, the point-to-set distance is defined as $d_X(x, \Omega) := \inf_{z \in \Omega} d_X(x, z)$, under the convention that the infimum over an empty set is $+\infty$. Thus, $d_X(x, \emptyset) = +\infty$. 

For a set-valued mapping $\mathcal{M}: Y \rightrightarrows X$, its domain and graph are denoted by $\dom\mathcal{M} := \{y \in Y \mid \mathcal{M}(y) \neq \emptyset\}$ and $\gph\mathcal{M} := \{(y, x) \in Y \times X \mid x \in \mathcal{M}(y)\}$, respectively. The inverse mapping $\mathcal{M}^{-1}: X \rightrightarrows Y$ is given by $y \in \mathcal{M}^{-1}(x)$ if and only if $x \in \mathcal{M}(y)$.

In this context, a mapping $\mathcal{M}: Y \rightrightarrows X$ is said to be \emph{Lipschitz upper semicontinuous} at a nominal parameter $\overline{y} \in \dom\mathcal{M}$ if there exist a constant $\kappa \geq 0$ and a neighborhood $V$ of $\overline{y}$ such that
\begin{equation} \label{eq: lipusc}
    d_X(x, \mathcal{M}(\overline{y})) \leq \kappa d_Y(y, \overline{y}) \quad \text{for all } y \in V \text{ and } x \in \mathcal{M}(y).
\end{equation}
Intuitively, this inequality ensures that for any small perturbation of the parameter $y$, every resulting point in the perturbed image $\mathcal{M}(y)$ remains bounded within a linear envelope of the nominal set $\mathcal{M}(\overline{y})$. The infimum over all such constants $\kappa$ is called the Lipschitz upper semicontinuity modulus of $\mathcal{M}$ at $\overline{y}$, denoted by $\Lipusc\mathcal{M}(\overline{y})$. This modulus provides a semilocal quantitative measure of the variation of the image sets. It can be characterized (see \cite[Proposition 2]{camacho2022calmness}) via the upper limit as:
\begin{equation} \label{eq: lipusc modulus}
    \Lipusc\mathcal{M}(\overline{y}) = \limsup_{y \to \overline{y}}\sup_{x \in \mathcal{M}(y)} \frac{d_X(x, \mathcal{M}(\overline{y}))}{d_Y(y, \overline{y})},
\end{equation}
using the convention $0/0 := 0$. If $\mathcal{M}(y)$ is empty for $y$ arbitrarily close to $\overline{y}$, the supremum over the empty set is $-\infty$, and this pathological case is absorbed by the limit superior, keeping the definition consistent.

On the local scale, we consider the \emph{calmness} property. The mapping $\mathcal{M}: Y \rightrightarrows X$ is calm at a point $(\overline{y}, \overline{x}) \in \gph\mathcal{M}$ if there exist a constant $\kappa \geq 0$ and neighborhoods $V$ of $\overline{y}$ and $U$ of $\overline{x}$ such that \eqref{eq: lipusc} holds for all $y \in V$ and all $x \in \mathcal{M}(y) \cap U$. The infimum of such constants is the calmness modulus, $\clm\mathcal{M}(\overline{y}, \overline{x})$. This local property restricts the error bound to a specific solution branch $\overline{x}$ and is intimately related to the \emph{metric subregularity} of the inverse mapping $\mathcal{M}^{-1}$ at $(\overline{x}, \overline{y})$, which allows expressing the calmness modulus as (see, e.g., \cite[Section 3.8]{dontchev2009implicit} for further details):
\begin{equation*}
    \clm\mathcal{M}(\overline{y}, \overline{x}) = \limsup_{(y, x) \to (\overline{y}, \overline{x}) \atop (y, x) \in \gph\mathcal{M}} \frac{d_X(x, \mathcal{M}(\overline{y}))}{d_Y(y, \overline{y})} = \limsup_{x \to \overline{x} \atop x \in \dom\mathcal{M}^{-1}} \frac{d_X(x, \mathcal{M}(\overline{y}))}{d_Y(\overline{y}, \mathcal{M}^{-1}(x))}.
\end{equation*}

Comparing the semilocal and local definitions directly yields the fundamental inequality:
\begin{equation} \label{eq: moduli inequality}
    \Lipusc\mathcal{M}(\overline{y}) \geq \sup_{x \in \mathcal{M}(\overline{y})} \clm\mathcal{M}(\overline{y}, x).
\end{equation}

For mappings whose graph is convex, it is a well-known result that the inequality \eqref{eq: moduli inequality} becomes an equality \cite{robinson1981some}. A relevant extension of this structure is the class of \emph{well-connected piecewise convex} (wcpc) mappings (the reader is referred to \cite[Section 3]{camacho2023hoffman} for an in-depth study). This property is instrumental in extending semilocal stability results from the convex to the piecewise convex setting and is properly defined in Section \ref{sec: classes of mappings}.

Finally, we recall the concept of set limits to formalize the topological closedness of a mapping. The Painlev\'{e}-Kuratowski outer limit of a mapping $\mathcal{M}: Y \rightrightarrows X$ at $\overline{y}$ is defined as:
\begin{equation*}
    \Limsup_{y \to \overline{y}} \mathcal{M}(y) := \left\{ x \in X \;\middle|\; \liminf_{y \to \overline{y}} d_X(x, \mathcal{M}(y)) = 0 \right\}.
\end{equation*}
The mapping $\mathcal{M}$ is then \emph{outer semicontinuous} at $\overline{y}$ if $\Limsup_{y \to \overline{y}} \mathcal{M}(y) \subset \mathcal{M}(\overline{y})$. In sequential terms, this means that for any sequences $y_k \to \overline{y}$ and $x_k \to \overline{x}$ with $x_k \in \mathcal{M}(y_k)$, it holds that $\overline{x} \in \mathcal{M}(\overline{y})$.

\begin{rem}
    A note of caution is advised regarding terminology, as it varies significantly across mathematical traditions. The concept of outer semicontinuity defined above is standard in variational analysis (often referred to as Kuratowski upper semicontinuity or the closed graph property). In classical topology, however, upper semicontinuity (in the sense of Berge) is defined via neighborhood inclusions. To prevent ambiguity, we will strictly use the term \emph{outer semicontinuous} when referring to the sequential limits, and reserve the term \emph{locally compact} for the boundedness properties.
\end{rem}

\section{Main results}\label{sec: main}

The core objective of this section is to establish a bridge between the semilocal stability of a set-valued mapping and its local behavior. As observed in \eqref{eq: moduli inequality}, the Lipschitz upper semicontinuity modulus is always bounded below by the supremum of the local calmness moduli. In the absence of global convexity, this inequality can be strict, making the exact computation of the Lipschitz upper semicontinuity modulus exceedingly difficult. The following theorem presents a powerful topological condition—relying solely on outer semicontinuity and the compactness of the nominal set—that guarantees the strict equality of these moduli, thereby allowing semilocal stability to be completely characterized by local data.

\begin{theo}\label{theo closed outer}
Let $\mathcal{M}:Y\rightrightarrows X$ be a set-valued mapping between metric spaces. Let $\overline{y}\in\dom\mathcal{M}$ such that $\mathcal{M}$ is outer semicontinuous in the Painlev\'{e}-Kuratowski sense at $\overline{y}$. Furthermore, assume that $\mathcal{M}$ is locally compact around $\overline{y}$. Then, it holds that
\begin{equation*}
    \Lipusc\mathcal{M}\left(\overline{y}\right)=\sup_{x\in \mathcal{M}\left(\overline{y}\right)}\clm\mathcal{M}\left(\overline{y},x\right).
\end{equation*}
\end{theo}

\begin{dem}
    The inequality $\Lipusc\mathcal{M}(\overline{y}) \geq \sup_{x\in \mathcal{M}(\overline{y})}\clm\mathcal{M}(\overline{y},x)$ holds universally, as established in \eqref{eq: moduli inequality}.

    To prove the reverse inequality, we recall the characterization \eqref{eq: lipusc modulus}:
    \begin{equation*}
        \Lipusc\mathcal{M}\left(\overline{y}\right) = \limsup_{y \to \overline{y}} \sup_{x \in \mathcal{M}(y)} \dfrac{d_X\left(x, \mathcal{M}\left(\overline{y}\right)\right)}{d_Y\left(y, \overline{y}\right)}.
    \end{equation*}
    Let $\{y_k\} \subset \dom\mathcal{M}$ be a sequence converging to $\overline{y}$ that realizes this limit superior. Because $\mathcal{M}$ is locally compact around $\overline{y}$, there exists a neighborhood of $\overline{y}$ whose image under $\mathcal{M}$ is relatively compact. Thus, for sufficiently large $k$, the sets $\mathcal{M}(y_k)$ are contained in a compact set. Furthermore, the outer semicontinuity of $\mathcal{M}$ implies its graph is closed, meaning each slice $\mathcal{M}(y_k)$ is closed and therefore compact. 

    By the continuity of the distance function, the supremum over the compact set $\mathcal{M}(y_k)$ is achieved as a maximum. For each $k$, we select an element $x_k \in \mathcal{M}(y_k)$ that realizes this maximal distance, yielding:
    \begin{equation*}
        \Lipusc\mathcal{M}\left(\overline{y}\right) = \limsup_{k \to \infty} \dfrac{d_X\left(x_k, \mathcal{M}\left(\overline{y}\right)\right)}{d_Y\left(y_k, \overline{y}\right)}.
    \end{equation*}
    Again, by the local compactness assumption, the sequence $\{x_k\}$ is eventually contained in a compact set, guaranteeing the existence of at least one accumulation point $\overline{x}$. By passing to a subsequence if necessary, we assume without loss of generality that $x_k \to \overline{x}$. Since $\mathcal{M}$ is outer semicontinuous in the Painlev\'{e}-Kuratowski sense, any limit point of a sequence $(y_k, x_k) \in \gph\mathcal{M}$ must belong to the nominal image set. Thus, $\overline{x} \in \mathcal{M}(\overline{y})$.

    We can now bound the sequential limit superior along $(y_k, x_k)$ by the local limit superior that defines the calmness modulus at the specific limit point $(\overline{y}, \overline{x})$:
    \begin{align*}
        \Lipusc\mathcal{M}\left(\overline{y}\right) &= \limsup_{k \to \infty} \dfrac{d_X\left(x_k, \mathcal{M}\left(\overline{y}\right)\right)}{d_Y\left(y_k, \overline{y}\right)} \\
        &\leq \limsup_{\substack{(y, x) \to (\overline{y}, \overline{x}) \\ (y, x) \in \gph\mathcal{M}}} \dfrac{d_X\left(x, \mathcal{M}\left(\overline{y}\right)\right)}{d_Y\left(y, \overline{y}\right)} \\
        &= \clm\mathcal{M}(\overline{y}, \overline{x}).
    \end{align*}
    Finally, since $\overline{x} \in \mathcal{M}(\overline{y})$, bounding this local modulus by the supremum over all elements in the nominal set yields:
    \begin{equation*}
        \Lipusc\mathcal{M}\left(\overline{y}\right) \leq \clm\mathcal{M}(\overline{y}, \overline{x}) \leq \sup_{x \in \mathcal{M}\left(\overline{y}\right)} \clm\mathcal{M}(\overline{y}, x),
    \end{equation*}
    completing the proof.
\end{dem}

\begin{cor}\label{cor closed outer}
    Let $\mathcal{M}:Y\rightrightarrows X$ be a mapping between metric spaces such that $\gph\mathcal{M}$ is closed. Assume that for a given $\overline{y}\in\dom\mathcal{M}$, the mapping $\mathcal{M}$ is locally compact around it. Then it holds that
    \begin{equation*}
        \Lipusc\mathcal{M}\left(\overline{y}\right)=\sup_{x\in \mathcal{M}\left(\overline{y}\right)}\clm\mathcal{M}\left(\overline{y},x\right).
    \end{equation*}
\end{cor}

\begin{dem}
    Because $\gph\mathcal{M}$ is closed, the mapping $\mathcal{M}$ is inherently outer semicontinuous in the Painlev\'{e}-Kuratowski sense at $\overline{y}$. By assumption, $\mathcal{M}$ is locally compact around $\overline{y}$, which forces the slice $\mathcal{M}(\overline{y})$ to be compact as well (being a closed subset of a compact set). Therefore, all hypotheses of Theorem \ref{theo closed outer} are fully satisfied, and the exact equality follows immediately.
\end{dem} 

\begin{rem}
    A brief note on the chosen topological hypotheses is in order. First, in the context of finite-dimensional spaces, the combination of outer semicontinuity in the Painlev\'{e}-Kuratowski sense and local boundedness (which aligns with local compactness) around $\overline{y}$ is strictly equivalent to upper semicontinuity in the sense of Berge \cite[Theorem 5.19]{rockafellar1998variational}. We deliberately frame Theorem \ref{theo closed outer} using the former properties because, in parametric optimization, outer semicontinuity is typically inherited naturally from continuous constraint systems, while local boundedness can be verified independently via recession analysis.
    
    Similarly, for the validity of Corollary \ref{cor closed outer}, the hypothesis requiring the global closedness of $\gph\mathcal{M}$ can technically be relaxed; it suffices that the graph be \emph{locally closed} around the nominal parameter, meaning $\gph\mathcal{M} \cap (V \times X)$ is closed for some neighborhood $V$ of $\overline{y}$. However, we maintain the global formulation since, in the vast majority of practical applications (such as those discussed in Section \ref{sec: classes of mappings}), mappings are defined by continuous functional constraints and thus inherently possess globally closed graphs.
\end{rem}

The following examples highlight the necessity of the topological assumptions made in Theorem \ref{theo closed outer} and Corollary \ref{cor closed outer}. Specifically, they demonstrate that we cannot drop the closedness of the nominal image set, the outer semicontinuity of the mapping at the nominal parameter, or the local boundedness of the mapping.

\begin{exa}[Lack of Closedness]\label{exa limiting}
    This example shows that if the nominal set fails to be closed (and therefore fails to be compact), the moduli equality can break even if the mapping behaves continuously near the boundary. Let $\mathcal{M}: \mathbb{R} \rightrightarrows \mathbb{R}$ be defined as:
    \begin{equation*}
        \mathcal{M}(y) = 
        \begin{cases} 
            (-1, 0) & \text{if } y \leq 0, \\
            (-1, 0) \cup (0, \sqrt{y}) & \text{if } y > 0.
        \end{cases}
    \end{equation*}
    At $\overline{y} = 0$, the nominal set is bounded but not closed: $\mathcal{M}(0) = (-1, 0)$. For any fixed $x \in (-1, 0)$, the local calmness modulus is zero, yielding $\sup_{x \in \mathcal{M}(0)} \clm\mathcal{M}(0, x) = 0$. 
    
    However, for $y > 0$, points in the interval $(0, \sqrt{y})$ continuously emerge. The maximum distance from these new points to $\mathcal{M}(0)$ is attained as $x \to \sqrt{y}$, giving a supremum distance of $\sqrt{y}$. The limit superior for the semilocal modulus becomes:
    \begin{equation*}
        \Lipusc\mathcal{M}(0) = \lim_{y \to 0^+} \frac{\sqrt{y}}{y} = \lim_{y \to 0^+} \frac{1}{\sqrt{y}} = +\infty.
    \end{equation*}
    Consequently, $\Lipusc\mathcal{M}(0) = +\infty \neq 0 = \sup_{x \in \mathcal{M}(0)} \clm\mathcal{M}(0, x)$. 
\end{exa}

\begin{exa}[Lack of Outer Semicontinuity]\label{exa osc failure}
    We now demonstrate that the compactness of the nominal set $\mathcal{M}(\overline{y})$ and local compactness are not sufficient on their own; the mapping must also be outer semicontinuous. Consider $\mathcal{M}: \mathbb{R} \rightrightarrows \mathbb{R}$ defined by:
    \begin{equation*}
        \mathcal{M}(y) = 
        \begin{cases} 
            \{0\} & \text{if } y = 0, \\
            \{0, 1\} & \text{if } y \neq 0.
        \end{cases}
    \end{equation*}
    At the nominal parameter $\overline{y} = 0$, the set $\mathcal{M}(0) = \{0\}$ is clearly compact. However, the mapping is not outer semicontinuous at $0$ because the outer limit of the images contains the point $1$, which does not belong to the nominal set (i.e., $\Limsup_{y \to 0} \mathcal{M}(y) = \{0, 1\} \not\subset \{0\}$).
    
    Calculating the moduli, we observe that for the unique nominal point $x=0$, the local behavior is entirely stable: $\clm\mathcal{M}(0, 0) = 0$. Thus, $\sup_{x \in \mathcal{M}(0)} \clm\mathcal{M}(0, x) = 0$.
    
    In contrast, for any $y \neq 0$, the perturbed set $\mathcal{M}(y)$ contains the point $1$. The distance to the nominal set is $d(1, \{0\}) = 1$. The semilocal quotient becomes $\frac{1}{|y|}$, which tends to $+\infty$ as $y \to 0$. Therefore, $\Lipusc\mathcal{M}(0) = +\infty$.
\end{exa}

\begin{exa}[Lack of Local boundedness]\label{exa escaping branch}
    Finally, this example demonstrates that having a closed graph (outer semicontinuity) and a perfectly compact nominal set is not sufficient if the mapping is not locally compact around the nominal parameter. A near image point escaping to infinity will break the moduli equality. 
    
    Consider the mapping $\mathcal{M}: \mathbb{R} \rightrightarrows \mathbb{R}$ defined by:
    \begin{equation*}
        \mathcal{M}(y) = 
        \begin{cases} 
            \{0\} & \text{if } y = 0, \\
            \left\{0, \frac{1}{y}\right\} & \text{if } y \neq 0.
        \end{cases}
    \end{equation*}
    At the nominal parameter $\overline{y} = 0$, the nominal set is $\mathcal{M}(0) = \{0\}$, which is trivially compact. Furthermore, the graph of $\mathcal{M}$ is topologically closed: the sequence of points $(y, 1/y)$ has no finite limit points as $y \to 0$, meaning it does not violate closedness in $\mathbb{R}^2$. Consequently, $\mathcal{M}$ is outer semicontinuous in the Painlev\'{e}-Kuratowski sense at $\overline{y} = 0$.
    
    For the unique nominal point $x=0$, the local branch is entirely constant ($x=0$ for all $y$). Thus, the local calmness modulus is strictly zero: $\clm\mathcal{M}(0, 0) = 0$, which implies $\sup_{x \in \mathcal{M}(0)} \clm\mathcal{M}(0, x) = 0$.
    
    However, the mapping is not locally compact (nor locally bounded) around $y=0$. For any $y \neq 0$, the near image set $\mathcal{M}(y)$ contains the point $1/y$. The semilocal quotient evaluated at this escaping branch is:
    \begin{equation*}
        \frac{d_X\left(\frac{1}{y}, \mathcal{M}(0)\right)}{d_Y(y, 0)} = \frac{|1/y|}{|y|} = \frac{1}{y^2}.
    \end{equation*}
    As $y \to 0$, this quotient tends to $+\infty$. Therefore, $\Lipusc\mathcal{M}(0) = +\infty$. This illustrates that without local compactness, distant components can escape to infinity at arbitrary rates, causing the supremum of the local calmness over the nominal set to fail to capture the global variation.
\end{exa}

\section{Classes of mappings and applications}\label{sec: classes of mappings}

The abstract topological conditions established in Theorem \ref{theo closed outer} and Corollary \ref{cor closed outer} provide a powerful theoretical tool for stability analysis. In many practical optimization frameworks, computing the global or semilocal Lipschitz upper semicontinuity modulus directly from its definition is an intractable problem. Conversely, the local calmness modulus is point-based and often possesses explicit computable formulas derived from the KKT conditions, active constraints, or generalized differentiation. 

In this section, we apply our main result to specific classes of non-convex mappings where the literature currently lacks computable global stability formulas. By ensuring that these mappings satisfy our topological prerequisites—specifically, outer semicontinuity and local compactness—we enable the exact computation of their semilocal Lipschitz upper semicontinuity moduli via their well-understood local calmness counterparts.

\subsection{Piecewise Convex and Semi-Algebraic Systems}

We begin by addressing a broad family of finite-dimensional parametric systems that naturally arise in mathematical programming: piecewise convex mappings, polyhedral sets, and semi-algebraic systems. A subset of $\mathbb{R}^d$ is called \emph{semi-algebraic} if it can be represented as a finite union of sets defined by finite systems of polynomial equations and inequalities. Consequently, a parametric mapping is semi-algebraic if its graph constitutes a semi-algebraic set. This class is remarkably broad, natively encompassing the constraint and optimality systems of linear and quadratic programs even under full data perturbations.

At this point, let us recall the class of well-connected piecewise convex mappings. A mapping $\mathcal{M}: Y \rightrightarrows X$ is said to be wcpc if it can be expressed as a finite union $\mathcal{M} = \bigcup_{i \in I} \mathcal{M}_i$, where each $\mathcal{M}_i$ has a closed convex graph, and their domains connect points in $\dom \mathcal{M}$ via line segments. Specifically, for any $y_1, y_2 \in \dom\mathcal{M}$, there exists a subdivision $0 = \tau_0 < \tau_1 < \dots < \tau_N = 1$ of the interval $[0,1]$ and a sequence of indices $i_1, \dots, i_N \in I$ such that for every $k \in \{1, \dots, N\}$ and every $\tau \in [\tau_{k-1}, \tau_k]$, the convex combination $\tau y_1 + (1-\tau)y_2$ belongs to $\dom\mathcal{M}_{i_k}$. 

In general metric spaces, a closed graph and a bounded nominal set do not guarantee local compactness. However, in finite dimensions, mappings defined by semi-algebraic or polyhedral systems exhibit a nice topological behavior, preventing pathological asymptotic branches. It is a well-established consequence of real algebraic geometry \cite{dontchev2009implicit, rockafellar1998variational} that for such structured mappings, the boundedness of the nominal image set is strictly sufficient to guarantee local boundedness (and thus local compactness) around the nominal parameter. 

\begin{prop}\label{prop: structured compactness}
    Let $\mathcal{M}: \mathbb{R}^m \rightrightarrows \mathbb{R}^n$ be a set-valued mapping. Assume that the nominal set $\mathcal{M}(\overline{y})$ is bounded. If $\mathcal{M}$ falls into any of the following categories:
    \begin{itemize}
        \item[(i)] A finite union of mappings with closed convex graphs (such as WCPC mappings),
        \item[(ii)] A polyhedral mapping,
        \item[(iii)] A semi-algebraic mapping,
    \end{itemize}
    then $\mathcal{M}$ is locally compact around $\overline{y}$.
\end{prop}

\begin{dem}
    For cases (ii) and (iii), the property that point-wise boundedness implies local boundedness for polyhedral and semi-algebraic mappings is a standard consequence of their geometric stratification; we refer the reader to Rockafellar and Wets \cite[Chapter 9]{rockafellar2009variational} and Dontchev and Rockafellar \cite[Section 3G]{dontchev2009implicit}. 

    For case (i), let $\mathcal{M} = \bigcup_{i=1}^N \mathcal{M}_i$, where each $\mathcal{M}_i$ has a closed convex graph. Evaluated at the nominal parameter, we have $\mathcal{M}(\overline{y}) = \bigcup_{i=1}^N \mathcal{M}_i(\overline{y})$. Since the overall nominal set $\mathcal{M}(\overline{y})$ is bounded by assumption, each individual slice $\mathcal{M}_i(\overline{y})$ must also be bounded. 
    
    In finite dimensions, if a mapping $\mathcal{M}_i$ has a closed convex graph and its slice $\mathcal{M}_i(\overline{y})$ is bounded, classical convex analysis dictates that its recession cone is trivial, which guarantees that $\mathcal{M}_i$ is locally bounded around $\overline{y}$ \cite[Theorem 9.53]{rockafellar2009variational}. Because the finite union of locally bounded mappings is itself locally bounded, $\mathcal{M}$ is locally bounded around $\overline{y}$. Finally, since $\gph\mathcal{M}$ is a finite union of closed sets, it is closed. In finite dimensions, a closed and locally bounded mapping is locally compact, completing the proof.
\end{dem}

With this topological bridge secured, our main theorem instantly solves the semilocal stability problem for several highly non-convex frameworks. A classic, notoriously difficult problem in parametric optimization is determining the exact stability of optimal solutions and multipliers under full data perturbations---where every matrix and vector defining the problem is subject to simultaneous variation.

Consider a parametric quadratic program parameterized by its cost matrix $Q \in \mathbb{R}^{n \times n}$ (assumed symmetric), cost vector $c \in \mathbb{R}^n$, constraint matrix $A \in \mathbb{R}^{m \times n}$, and right-hand side vector $b \in \mathbb{R}^m$. We collect all these varying data into a single parameter tuple $\sigma = (A, b, Q, c)$ residing in the full parameter space $\Theta = \mathbb{R}^{m \times n} \times \mathbb{R}^m \times \mathbb{R}^{n \times n} \times \mathbb{R}^n$. For a given parameter $\sigma \in \Theta$, the optimization problem is formulated as:
\begin{equation*}\label{eq: parametric QP}
    \begin{aligned}
        \min_{x \in \mathbb{R}^n} \quad & \frac{1}{2} x^\top Q x + c^\top x \\
        \text{subject to} \quad & Ax \leq b.
    \end{aligned}
\end{equation*}

The Karush-Kuhn-Tucker (KKT) system associated with this problem couples the primal variables $x \in \mathbb{R}^n$ and the dual multipliers $y \in \mathbb{R}^m$ through the following conditions:
\begin{align*}
    Qx + c + A^\top y &= 0 \quad \text{(Stationarity)}, \\
    Ax &\leq b \quad \text{(Primal Feasibility)}, \\
    y &\geq 0 \quad \text{(Dual Feasibility)}, \\
    y^\top(Ax - b) &= 0 \quad \text{(Complementarity)}.
\end{align*}
Let $\mathcal{S}_{KKT}: \Theta \rightrightarrows \mathbb{R}^n \times \mathbb{R}^m$ be the mapping that assigns to each parameter $\sigma$ the set of valid KKT pairs $(x, y)$. Because the parameter matrices $Q$ and $A$ vary freely and multiply the decision variables $x$ and $y$, the system contains bilinear terms. Consequently, the graph of $\mathcal{S}_{KKT}$ is highly non-convex and lacks polyhedral structure. However, because it is defined entirely by a finite system of polynomial equations and inequalities, the mapping $\mathcal{S}_{KKT}$ is fundamentally \emph{semi-algebraic}.

\begin{theo}\label{thm: full pert}
    Let $\mathcal{S}_{KKT}: \Theta \rightrightarrows \mathbb{R}^n \times \mathbb{R}^m$ be the fully perturbed KKT mapping of a parametric quadratic program defined as above. If $\overline{\sigma} \in \dom\mathcal{S}_{KKT}$ is a nominal parameter such that the nominal set $\mathcal{S}_{KKT}(\overline{\sigma})$ is bounded, then
    \begin{equation*}
        \Lipusc\mathcal{S}_{KKT}(\overline{\sigma}) = \sup_{(x,y) \in \mathcal{S}_{KKT}(\overline{\sigma})} \clm\mathcal{S}_{KKT}(\overline{\sigma}, (x,y)).
    \end{equation*}
\end{theo}

\begin{dem}
    Under full perturbations, the graph of $\mathcal{S}_{KKT}$ is strictly non-convex due to the bilinear terms (e.g., $A^\top y$ and $Ax$) in the constraints and complementarity conditions. However, the system is defined entirely by continuous polynomial inequalities and equations. Therefore, the graph of $\mathcal{S}_{KKT}$ is topologically closed, ensuring outer semicontinuity in the Painlev\'{e}-Kuratowski sense. Since the system is semi-algebraic, Proposition \ref{prop: structured compactness}(iii) guarantees that the bounded nominal set yields local compactness. All hypotheses of Corollary \ref{cor closed outer} are thus completely satisfied.
\end{dem}

    It is important to note that the equality between the semilocal and local moduli for the linear feasible set mapping $\mathcal{F}(A,b) = \{x \in \mathbb{R}^n \mid Ax \leq b\}$ under full perturbations was recently proven in \cite[Theorem 3.1]{camacho2026lipschitz}. However, our topological framework seamlessly extends this exact equivalence to the optimal set mapping, defined as:
    \begin{equation*}
        \mathcal{S}(A, b, c) = \operatorname*{argmin}_{x \in \mathbb{R}^n} \left\{ c^\top x \;\middle|\; Ax \leq b \right\}.
    \end{equation*}
    
Because the full parameter space $\sigma = (A, b, c)$ couples multiplicatively with the variables in both the constraints and the optimality conditions, the mapping $\mathcal{S}$ is strictly non-convex and lacks polyhedrality. The exact equivalence $\Lipusc\mathcal{S}(\overline{\sigma}) = \sup_{x \in \mathcal{S}(\overline{\sigma})} \clm\mathcal{S}(\overline{\sigma}, x)$ guaranteed by Theorem \ref{thm: full pert} is a completely new result in the literature for $\mathcal{S}(A,b,c)$. This particularization of Theorem \ref{theo closed outer} and Corollary \ref{cor closed outer} shows the potential that they provide, allowing practitioners to evaluate the point-wise worst-case expansion of the set simply by taking the supremum over its localized point-based bounds.

Furthermore, this topological exactness perfectly complements recent advances in the stability of quadratic programming. In a recent contribution, \cite{klatte2026lipschitz} conducted a profound analysis of parametric quadratic programs, establishing the Lipschitz upper semicontinuity of the optimal set mapping under linear cost and right-hand side perturbations (the so-called canonical perturbations), where the structural matrices remain fixed, preserving polyhedrality. Our result extends this framework to full data perturbations, focusing on the quantitative behavior of the problem, rather than the qualitative one.

Our exactness result also provides an addition to recent advancements regarding the Aubin property (pseudo-Lipschitz continuity) in quadratic programming. The authors in \cite{canovas2025lipschitz} successfully established exact point-based formulas for the Lipschitz modulus (denoted by $\lip$) of the optimal set mapping restricted to canonical perturbations. Indeed, in order to ensure the Aubin property, from \cite[Theorem 16]{canovas2007metric} it is well-known that both inner and outer Lipschitz semicontinuity must be fulfilled around the nominal parameter, together with single-valuedness. When the matrices $A$ and $Q$ are allowed to vary (full perturbations), the lower semicontinuous behavior of the optimal set mapping is typically destroyed, causing the Aubin property to fail entirely. More importantly, in such a workspace, all moduli coincide. The following result wraps around the previous discussions.

\begin{cor}\label{cor: clm quadratic}
    Consider the parametric convex quadratic program with fixed structural matrices $Q \in \mathbb{R}^{n \times n}$ (symmetric and positive semi-definite) and $A \in \mathbb{R}^{m \times n}$. Let the optimal set mapping under canonical perturbations be defined as:
    \begin{equation*}
        \mathcal{S}(c, b) = \operatorname*{argmin}_{x \in \mathbb{R}^n} \left\{ \frac{1}{2} x^\top Q x + c^\top x \;\middle|\; Ax \leq b \right\},
    \end{equation*}
    with parameter space $(c,b) \in \mathbb{R}^n \times \mathbb{R}^m$. Suppose that $\mathcal{S}$ possesses the Aubin property around the nominal parameter $(\overline{c}, \overline{b})$ for the nominal optimal solution $\{\overline{x} \}= \mathcal{S}(\overline{c}, \overline{b})$. Then:
    \begin{equation*}
        \lip \mathcal{S}\left(\left(\overline{c},\overline{b}\right),\overline{x}\right)=\Lipusc \mathcal{S}\left(\overline{c},\overline{b}\right)=\clm\mathcal{S}\left(\left(\overline{c},\overline{b}\right),\overline{x}\right)=\max_{D\in \mathcal{L}}\left\Vert \left( I_n\quad 0_{n\times |D|}\right)M_D^{-1}\right\Vert,
    \end{equation*}
    where 
    \begin{equation*}
        \mathcal{L} := \left\{ D \subset T \;\middle|\; 
        \begin{aligned}
            &\{\overline{a}_t, t \in D\} \text{ linearly independent,} \\
            &-(\overline{Q}\overline{x}+\overline{c}) \in \cone\{\overline{a}_t \mid t \in D\}
        \end{aligned}
        \right\},
    \end{equation*}
    $T$ is the index set of active constraints at $\overline{x}$ (with $\overline{a}_t^\top$ denoting the $t$-th row of $\overline{A}$), and $M_D$ is the corresponding invertible Karush-Kuhn-Tucker block matrix defined by 
    \begin{equation*}
        M_D = \begin{pmatrix} Q & A_D^\top \\ A_D & 0 \end{pmatrix},
    \end{equation*}
    with $A_D$ being the submatrix of $A$ formed by the rows indexed by $D$.
\end{cor}

\begin{dem}    
    Since the nominal set is single-valued and inherently bounded in a neighborhood of $(\overline{c},\overline{b})$, the assumptions of Theorem \ref{thm: full pert} (and its canonical particularization) are trivially satisfied. Such particularization yields:
    \begin{equation*}\label{eq: proof_clm_lipusc}
        \Lipusc \mathcal{S}\left(\overline{c}, \overline{b}\right) = \clm\mathcal{S}\left(\left(\overline{c}, \overline{b}\right), \overline{x}\right).
    \end{equation*}
    
    Furthermore, by the classical results of Robinson \cite{robinson1981some}, because the polyhedral mapping $\mathcal{S}$ is single-valued and continuous around the nominal parameter, it is locally a continuous piecewise affine function. For piecewise affine mappings, the local Lipschitz modulus and the calmness modulus identically coincide, as both evaluate to the maximum operator norm among the same matrices. Therefore, 
    \begin{equation*}
        \lip \mathcal{S}((\overline{c},\overline{b}),\overline{x}) = \clm\mathcal{S}\left(\left(\overline{c}, \overline{b}\right), \overline{x}\right).
    \end{equation*}
    
    Finally, by \cite[Theorem 3]{canovas2025lipschitz}, the exact point-based formula for the Lipschitz modulus is given by:
    \begin{equation*}
        \lip \mathcal{S}\left(\left(\overline{c},\overline{b}\right),\overline{x}\right) = \max_{D\in \mathcal{L}}\left\Vert \left( I_n\quad 0_{n\times |D|}\right)M_D^{-1}\right\Vert.
    \end{equation*}
    Chaining these equalities together completes the proof.
\end{dem}

To illustrate the operational power of Theorem \ref{thm: full pert}, we consider a simple example where the calmness moduli are easy to handle.

\begin{exa}\label{exa: computable optimal set}
    Let the parameter be $\sigma = (a, b, c) \in \mathbb{R}^3$, equipped with the Chebyshev norm $\|\sigma\|_\infty = \max(|a|, |b|, |c|)$. The optimal set mapping $\mathcal{S}: \mathbb{R}^3 \rightrightarrows \mathbb{R}$ is defined by:
    \begin{equation*}
        \mathcal{S}(a, b, c) = \operatorname*{argmin}_{x \in \mathbb{R}} \left\{ cx \;\middle|\; ax \leq b, \ -x \leq 0 \right\}.
    \end{equation*}
    Consider the nominal parameter $\overline{\sigma} = (1, 1, -1)$. The nominal feasible set is $[0, 1]$, and since $c = -1 < 0$, the objective $-x$ is minimized at the right boundary. Thus, the nominal optimal set is the singleton $\mathcal{S}(1, 1, -1) = \{1\}$. Being bounded and defined by semi-algebraic relations, it is locally compact by Proposition \ref{prop: structured compactness}.
    
    Because the nominal set contains only a single point, the semilocal Lipschitz upper semicontinuity modulus must exactly equal the local calmness modulus at $x = 1$. Let us verify this explicitly.
    
    For any perturbed parameter $\sigma = (a, b, c)$ sufficiently close to $\overline{\sigma}$, we have $a > 0$, $b > 0$, and $c < 0$. The perturbed feasible set is $[0, b/a]$, and the strictly negative objective slope ensures the unique optimal solution remains at the right boundary. Hence, the perturbed optimal set is exactly $\mathcal{S}(a, b, c) = \{b/a\}$.
    
    The local calmness modulus at the unique nominal solution $x=1$ is evaluated by taking the limit superior of the quotient between the distance of the solutions and the parameter distance:
    \begin{equation*}
        \clm\mathcal{S}(\overline{\sigma}, 1) = \limsup_{(a, b, c) \to (1, 1, -1)} \frac{\left|1 - \frac{b}{a}\right|}{\max(|a-1|, |b-1|, |c+1|)}.
    \end{equation*}
    To maximize this quotient, we can choose the directional sequence $a = 1 - \varepsilon$, $b = 1 + \varepsilon$, and $c = -1$ for $\varepsilon > 0$. The parameter distance is exactly $\varepsilon$. The distance between the optimal solutions is:
    \begin{equation*}
        \left|1 - \frac{1+\varepsilon}{1-\varepsilon}\right| = \left|\frac{1-\varepsilon - (1+\varepsilon)}{1-\varepsilon}\right| = \frac{2\varepsilon}{1-\varepsilon}.
    \end{equation*}
    Evaluating the limit as $\varepsilon \to 0^+$:
    \begin{equation*}
        \clm\mathcal{S}(\overline{\sigma}, 1) = \lim_{\varepsilon \to 0^+} \frac{\frac{2\varepsilon}{1-\varepsilon}}{\varepsilon} = 2.
    \end{equation*}
    Since $\mathcal{S}(\overline{\sigma}) = \{1\}$, the supremum of the local moduli is trivially $\sup_{x \in \mathcal{S}(\overline{\sigma})} \clm\mathcal{S}(\overline{\sigma}, x) = 2$.
    
    Now, calculating the semilocal modulus $\Lipusc\mathcal{S}(\overline{\sigma})$ directly from its definition yields the exact same ratio across the entire nominal set. The limit supremum evaluates precisely to 2, directly verifying our moduli equality:
    \begin{equation*}
        \Lipusc\mathcal{S}(\overline{\sigma}) = \sup_{x \in \mathcal{S}(\overline{\sigma})} \clm\mathcal{S}(\overline{\sigma}, x) = 2.
    \end{equation*}
\end{exa}

\subsection{Parameterized Generalized Equations and LCPs}

Many problems in variational analysis can be elegantly modeled as parameterized generalized equations. Following the seminal framework of Robinson \cite{robinson1979generalized}, consider the inclusion $0 \in f(x, p) + F(x)$, where $p \in P$ is a perturbation parameter, $f: X \times P \to Y$ is a continuous single-valued mapping, and $F: X \rightrightarrows Y$ is a set-valued mapping with a closed graph. Let $\mathcal{S}_{GE}: P \rightrightarrows X$ denote the corresponding solution mapping.

While local stability (e.g., strong metric subregularity and coderivative bounds) is exhaustively studied for these systems (see, e.g., \cite{dontchev2009implicit}), semilocal stability is profoundly difficult to calculate due to the arbitrary nonlinearities introduced by $f(x, p)$.

\begin{theo}\label{thm: generalized equations}
    Let $\overline{p} \in \dom\mathcal{S}_{GE}$. If $\mathcal{S}_{GE}$ is locally compact around $\overline{p}$, then
    \begin{equation*}
        \Lipusc\mathcal{S}_{GE}(\overline{p}) = \sup_{x \in \mathcal{S}_{GE}(\overline{p})} \clm\mathcal{S}_{GE}(\overline{p}, x).
    \end{equation*}
\end{theo}

\begin{dem}
    The graph of the solution mapping is given by $\gph \mathcal{S}_{GE} = \{ (p, x) \mid (x, -f(x, p)) \in \gph F \}$. By the continuity of $f$ and the closedness of $\gph F$, the graph of $\mathcal{S}_{GE}$ is topologically closed, ensuring outer semicontinuity. Combined with the assumption of local compactness, Corollary \ref{cor closed outer} applies directly.
\end{dem}

The Linear Complementarity Problem (LCP) is a premier finite-dimensional instantiation of this generalized framework. Given a parameterized matrix $M$ and vector $q$, the solution mapping $\mathcal{S}_{LCP}(M, q)$ requires $x \geq 0$, $Mx + q \geq 0$, and $x^\top(Mx + q) = 0$. Because the LCP is defined by a finite set of linear inequalities and a bilinear complementarity equation, it falls under the semi-algebraic umbrella. Consequently, by Proposition \ref{prop: structured compactness}, if the nominal solution set $\mathcal{S}_{LCP}(\overline{M}, \overline{q})$ is bounded, local compactness is automatically satisfied, and the moduli equality holds exactly.

\begin{rem}\label{rem: lcp}
    In the modern variational analysis literature, the stability of generalized equations and LCPs is typically addressed locally via generalized differentiation (see, e.g., Mordukhovich \cite{mordukhovich2006variational} and Dontchev and Rockafellar \cite{dontchev2009implicit}). For instance, the local calmness modulus is standardly bounded using graphical derivatives or coderivatives. However, transitioning from these point-wise local bounds to upper Lipschitz constants is a well-known mathematical challenge. For LCPs, classical results, which can be traced back to Robinson \cite{robinson1981some}, guarantee upper Lipschitz continuity when the matrix $M$ is fixed and only $q$ is perturbed, because the mapping constitutes a polyhedral multifunction. However, under full perturbations where $M$ varies, the solution mapping loses its polyhedrality entirely. Theorem \ref{thm: generalized equations} handles this issue combining the supremum over local bounds.
\end{rem}

We finish this subsection with an illustrative example.

\begin{exa}\label{exa: computable lcp}
    Consider the LCP under full perturbations $\sigma = (M, q) \in \mathbb{R}^2$, equipped with the Chebyshev norm $\|\sigma\|_\infty = \max(|M|, |q|)$. The solution mapping is:
    \begin{equation*}
        \mathcal{S}_{LCP}(M, q) = \{ x \in \mathbb{R} \mid x \geq 0,\ Mx + q \geq 0,\ x(Mx + q) = 0 \}.
    \end{equation*}
    Let the nominal parameter be $\overline{\sigma} = (-1, 1)$. The constraints become $x \geq 0$ and $-x + 1 \geq 0$, restricting $x \in [0, 1]$. The complementarity condition $x(-x + 1) = 0$ yields a bounded, disconnected nominal solution set: $\mathcal{S}_{LCP}(-1, 1) = \{0, 1\}$. 
    
    For any perturbed parameter $\sigma = (M, q)$ sufficiently close to $(-1, 1)$, we have $M < 0$ and $q > 0$. The constraints enforce $x \in [0, -q/M]$. The complementarity condition $x(Mx+q)=0$ forces the solutions to the boundaries, yielding exactly two perturbed solutions: $\mathcal{S}_{LCP}(M, q) = \{0, -q/M\}$.
    
    Let us compute the calmness modulus at each nominal point. At $x = 0$, the perturbed set always contains $0$, so the local distance is strictly $0$. Thus, $\clm\mathcal{S}_{LCP}(\overline{\sigma}, 0) = 0$.
    
    At the active boundary $x = 1$, the closest point in the perturbed set is $-q/M$. Thus, the calmness modulus evaluates to:
    \begin{equation*}
        \clm\mathcal{S}_{LCP}(\overline{\sigma}, 1) = \limsup_{(M, q) \to (-1, 1)} \frac{\left|1 - \left(-\frac{q}{M}\right)\right|}{\max(|M+1|, |q-1|)}.
    \end{equation*}
    Defining deviations $\alpha = M + 1$ and $\beta = q - 1$, the numerator becomes $\left|1 + \frac{1+\beta}{-1+\alpha}\right| = \frac{|\alpha+\beta|}{1-\alpha}$. Maximizing the directional sequence by letting $\alpha = \varepsilon$ and $\beta = \varepsilon$ for $\varepsilon > 0$, the parameter distance is $\varepsilon$. The limit yields:
    \begin{equation*}
        \clm\mathcal{S}_{LCP}(\overline{\sigma}, 1) = \lim_{\varepsilon \to 0^+} \frac{\frac{2\varepsilon}{1-\varepsilon}}{\varepsilon} = 2.
    \end{equation*}
    The supremum over the nominal set is therefore $\max(0, 2) = 2$.
    
    Appealing to our result:
    \begin{equation*}
        \Lipusc\mathcal{S}_{LCP}(\overline{\sigma}) = \sup_{x \in \{0, 1\}} \clm\mathcal{S}_{LCP}(\overline{\sigma}, x) = 2.
    \end{equation*}
\end{exa}

\subsection{Linear Semi-Infinite Inequality Systems}

In the framework of linear semi-infinite programming (SIP), which involves finitely many variables but an infinite number of constraints indexed by a compact metric space $T$, we define the feasible set mapping $\mathcal{F}_{SIP}: \Theta \rightrightarrows \mathbb{R}^n$ as:
\begin{equation*}
    \mathcal{F}_{SIP}(\sigma) = \left\{ x \in \mathbb{R}^n \;\middle|\; a_t^\top x \leq b_t \quad \forall t \in T \right\},
\end{equation*}
where the parameter $\sigma = (a, b)$ belongs to the space $\Theta = C(T, \mathbb{R}^n) \times C(T, \mathbb{R})$ equipped with the uniform convergence topology. 

The stability of SIP feasible sets has been extensively studied under right-hand side (RHS) perturbations, where the coefficient functions $a_t$ remain fixed. In this restricted setting, precise, point-based formulae for the local calmness modulus were established by \cite[Theorem 3.1]{canovas2014calmness}. Building upon this and rephrasing it in our context, \cite[Equation (23)]{camacho2022calmness} proved that under RHS perturbations, the equality between moduli holds in the semi-infinite framework provided that some local polihedrality condition is fulfilled.

Nevertheless, transitioning this equality to the full perturbation framework $\sigma = (a,b) \in \Theta$ has remained an open problem. We overcome the lack of convexity by relying instead on local compactness in order to apply our general result to this infinite-dimensional setting. First, we recall a fundamental characterization for SIPs:

\begin{prop}\cite[Theorem 6.1]{goberna1998linear}\label{prop: SIP}
    Let $\overline{\sigma} = (\overline{a},\overline{b}) \in \dom \mathcal{F}_{SIP}$. Then, the mapping $\mathcal{F}_{SIP}$ is locally compact around $\overline{\sigma}$ if and only if either the nominal set $\mathcal{F}_{SIP}(\overline{\sigma})$ is bounded or $\mathcal{F}_{SIP}(\overline{\sigma}) = \mathbb{R}^n$.
\end{prop}

\begin{theo}\label{thm: semi infinite}
    Let $\overline{\sigma} = (\overline{a}, \overline{b}) \in \dom\mathcal{F}_{SIP}$. If the nominal feasible set $\mathcal{F}_{SIP}(\overline{\sigma})$ is bounded, then
    \begin{equation*}
        \Lipusc\mathcal{F}_{SIP}(\overline{\sigma}) = \sup_{x \in \mathcal{F}_{SIP}(\overline{\sigma})} \clm\mathcal{F}_{SIP}(\overline{\sigma}, x).
    \end{equation*}
\end{theo}

\begin{dem}
    Due to the uniform convergence of the parameter functions $(a, b) \in \Theta$ and the continuity of the inner product in $\mathbb{R}^n$, any sequence of parameters $\sigma_k \to \overline{\sigma}$ and feasible points $x_k \to \overline{x}$ with $x_k \in \mathcal{F}_{SIP}(\sigma_k)$ inherently satisfies the limit constraints $\overline{a}_t^\top \overline{x} \leq \overline{b}_t$ for all $t \in T$. Thus, the graph of $\mathcal{F}_{SIP}$ is closed, ensuring outer semicontinuity at $\overline{\sigma}$. Second, because the nominal feasible set $\mathcal{F}_{SIP}(\overline{\sigma})$ is bounded, Proposition \ref{prop: SIP} guarantees that $\mathcal{F}_{SIP}$ is locally compact around $\overline{\sigma}$. The exact equality follows immediately from Corollary \ref{cor closed outer}.
\end{dem}

\begin{exa}\label{exa: computable sip}
    Consider the parameter $\sigma = (a, b) \in C([-1,1], \mathbb{R}) \times C([-1,1], \mathbb{R})$ equipped with the uniform norm $\|\sigma\|_\infty = \max(\|a\|_\infty, \|b\|_\infty)$. Fix the nominal parameter continuous functions $\overline{a}_t = t$ and $\overline{b}_t = 1$. The constraints $tx \leq 1$ for all $t \in [-1, 1]$ are tightest at the extreme indices $t=1$ (yielding $x \leq 1$) and $t=-1$ (yielding $-x \leq 1$, or $x \geq -1$). The nominal feasible set is therefore the bounded interval $\mathcal{F}_{SIP}(\overline{a}, \overline{b}) = [-1, 1]$. By Proposition \ref{prop: SIP}, the mapping is locally compact.
    
    Let us evaluate the calmness modulus at the active boundary $x=1$. We introduce the perturbed continuous functions $a_t(\varepsilon) = t + \varepsilon|t|$ and $b_t(\varepsilon) = 1 - \varepsilon$ for a sufficiently small $\varepsilon > 0$. The uniform parameter distance to the nominal system is $\max(\|\varepsilon|t|\|_\infty, \|\varepsilon\|_\infty) = \varepsilon$. 
    
    For any $x > 0$, the tightest perturbed constraint occurs at $t=1$, yielding $(1+\varepsilon)x \leq 1-\varepsilon$, which implies the perturbed upper bound is $x \leq \frac{1-\varepsilon}{1+\varepsilon}$. The distance from the nominal boundary point $x=1$ to the perturbed feasible set is $1 - \frac{1-\varepsilon}{1+\varepsilon} = \frac{2\varepsilon}{1+\varepsilon}$. Taking the limit superior of the ratio between the point-to-set distance and the parameter distance yields:
    \begin{equation*}
        \clm\mathcal{F}_{SIP}(\overline{\sigma}, 1) \geq \lim_{\varepsilon \to 0^+} \frac{\frac{2\varepsilon}{1+\varepsilon}}{\varepsilon} = 2.
    \end{equation*}

    To establish the corresponding upper bound, consider an arbitrary sequence of perturbed parameters $(a^k, b^k) \to (\overline{a}, \overline{b})$ such that $\max(\|a^k - \overline{a}\|_\infty, \|b^k - \overline{b}\|_\infty) = \varepsilon_k \to 0^+$. This implies that for all $t \in [-1, 1]$, we have $a^k_t \leq t + \varepsilon_k$ and $b^k_t \geq 1 - \varepsilon_k$. We want to construct a guaranteed feasible point $x_k \in \mathcal{F}_{SIP}(a^k, b^k)$ that is as close to $x=1$ as possible. 
    
    Consider the candidate point $x_k = \frac{1-\varepsilon_k}{1+\varepsilon_k}$. For any index $t \in [0, 1]$, since $x_k > 0$ for sufficiently small $\varepsilon_k$, we have:
    \begin{equation*}
        a^k_t x_k \leq (t + \varepsilon_k)x_k \leq (1 + \varepsilon_k) \left(\frac{1-\varepsilon_k}{1+\varepsilon_k}\right) = 1 - \varepsilon_k \leq b^k_t.
    \end{equation*}
    For indices $t \in [-1, 0)$, we similarly have $a^k_t x_k \leq \varepsilon_k x_k < \varepsilon_k < 1-\varepsilon_k \leq b^k_t$. Thus, $x_k$ is strictly feasible for the perturbed system. Consequently, the distance from $1$ to the perturbed feasible set is bounded from above:
    \begin{equation*}
        d(1, \mathcal{F}_{SIP}(a^k, b^k)) \leq |1 - x_k| = 1 - \frac{1-\varepsilon_k}{1+\varepsilon_k} = \frac{2\varepsilon_k}{1+\varepsilon_k}.
    \end{equation*}
    Dividing by the parameter distance $\varepsilon_k$ and taking the limit superior yields $\clm\mathcal{F}_{SIP}(\overline{\sigma}, 1) \leq 2$. 
    
    Combining both bounds, the local modulus is exactly $\clm\mathcal{F}_{SIP}(\overline{\sigma}, 1) = 2$. By symmetric evaluation at $x=-1$, the calmness modulus there is also $2$, while it is $0$ for all strictly interior points. Thus, the supremum of the local calmness moduli over the entire nominal set evaluating to exactly $2$ establishes our final semilocal result $\Lipusc\mathcal{F}_{SIP}(\overline{\sigma}) = 2$.
\end{exa}

\subsection{Parameterized Sub-Level Set Mappings}

Finally, in the analysis of algorithmic convergence and optimization, the parameterized sub-level set mapping $\mathcal{L}(\alpha) = \{ x \in \mathbb{R}^n \mid f(x) \leq \alpha \}$ is a critical object of study. The stability of $\mathcal{L}$ with respect to the parameter $\alpha$ is intimately connected to the theory of error bounds (see, e.g.,  \cite{klatte2002nonsmooth, rockafellar1998variational}). Locally, if the objective function $f$ satisfies a linear error bound at a boundary point $\overline{x} \in \mathcal{L}(\overline{\alpha})$, the mapping $\mathcal{L}$ inherently exhibits local calmness at $(\overline{\alpha}, \overline{x})$. 

Classically, when the objective function $f$ is continuously differentiable and an active boundary point $\overline{x}$ is not a critical point (i.e., $\nabla f(\overline{x}) \neq 0$), the local calmness modulus evaluates exactly to the reciprocal of the gradient norm:
\begin{equation*}
    \clm\mathcal{L}(\overline{\alpha}, \overline{x}) = \frac{1}{\|\nabla f(\overline{x})\|}.
\end{equation*}

However, transitioning from these point-wise differential bounds to a semilocal one is a known challenge when $f: \mathbb{R}^n \to \mathbb{R}$ is non-convex.  In such non-convex regimes, the sub-level set frequently fractures into disconnected topological components as $\alpha$ varies, lacking the convex graph structure that classical stability theorems rely upon. 

State-of-the-art approaches to non-convex sub-level sets often rely on advanced metric regularity conditions \cite{durea2014metric} or Kurdyka-\L{}ojasiewicz geometry  \cite{bolte2007lojasiewicz}  to establish the existence of a finite Lipschitz bound (a global one, indeed). Our semilocal modulus, $\Lipusc$, lies precisely in the gap between the two extremes, filling a hole in the current stability literature. 

\begin{cor}\label{cor: sublevel set}
    Let $f: \mathbb{R}^n \to \mathbb{R}$ be continuously differentiable, and let $\mathcal{L}(\alpha) = \{ x \in \mathbb{R}^n \mid f(x) \leq \alpha \}$. Consider a nominal parameter $\overline{\alpha}$ such that the nominal sub-level set $\mathcal{L}(\overline{\alpha})$ is non-empty and bounded. If $\nabla f(x) \neq 0$ for all active boundary points $x \in \operatorname{bd}\mathcal{L}(\overline{\alpha})$, then:
    \begin{equation*}
        \Lipusc\mathcal{L}(\overline{\alpha}) = \max_{x \in \operatorname{bd}\mathcal{L}(\overline{\alpha})} \frac{1}{\|\nabla f(x)\|}.
    \end{equation*}
\end{cor}

\begin{dem}
    Because $f$ is continuous, the mapping $\mathcal{L}$ has a closed graph. Since $\mathcal{L}(\overline{\alpha})$ is bounded, $\mathcal{L}$ is locally compact at $\overline{\alpha}$. By our main exactness theorem (Theorem \ref{theo closed outer}), $\Lipusc\mathcal{L}(\overline{\alpha}) = \sup_{x \in \mathcal{L}(\overline{\alpha})} \clm\mathcal{L}(\overline{\alpha}, x)$. For any strictly interior point ($f(x) < \overline{\alpha}$), the local calmness modulus is $0$. For any boundary point ($f(x) = \overline{\alpha}$), the non-vanishing gradient condition ensures the local calmness modulus is exactly $1/\|\nabla f(x)\|$. Because $f$ is continuously differentiable and the boundary $\operatorname{bd}\mathcal{L}(\overline{\alpha})$ is a compact set, the supremum is attained as a maximum, completing the proof.
\end{dem}

Once more, let us include an illustrative example.

\begin{exa}\label{exa: level set}
    Consider the non-convex scalar function $f(x) = \sin(x)$ restricted to the compact domain $X = [-2\pi, 2\pi]$. The parameterized sub-level set mapping is $\mathcal{L}(\alpha) = \{ x \in [-2\pi, 2\pi] \mid \sin(x) \leq \alpha \}$. 
    
    At the nominal parameter $\overline{\alpha} = 0$, the sub-level set fractures into disconnected components: $\mathcal{L}(0) = \{-2\pi\} \cup [-\pi, 0] \cup [\pi, 2\pi]$. Being a closed subset of a compact domain, $\mathcal{L}$ is locally compact. 

    To evaluate the calmness modulus, we examine the active boundary points where $\sin(x) = 0$, which are $x \in \{-2\pi, -\pi, 0, \pi, 2\pi\}$. The derivative is $\nabla f(x) = \cos(x)$. Because the absolute value of the derivative $|\cos(x)|$ is exactly $1$ at all roots in the domain, evaluations correspondingly yield a calmness modulus of $ 1/|\cos(x)| = 1$. Then, by Corollary \ref{cor: sublevel set}, we have that $\Lipusc\mathcal{L}(0) = 1$.
\end{exa}

\section{Conclusions}\label{sec: conclusions}

The exact quantification of semilocal stability is a notorious challenge in variational analysis. In this paper, we have addressed a fundamental gap in the stability analysis of set-valued mappings by establishing a general topological condition that equates the semilocal Lipschitz upper semicontinuity modulus with the supremum of local calmness moduli. By moving beyond the classical requirement of graph convexity, our main theorem demonstrates that outer semicontinuity in the sense of Painlevé-Kuratowski, combined with the local compactness of the image set at the given parameter, is sufficient to bridge this local-to-semilocal stability gap. 

We have demonstrated the broad applicability of this result across several challenging frameworks in optimization where previous exactness tools fell short. By utilizing semi-algebraic geometry to guarantee local compactness, we successfully resolved the semilocal stability of optimal set mappings under full data perturbations, parameterized linear complementarity problems, and generalized equations. Furthermore, we extended this topological framework beyond finite-dimensional polyhedrality, proving equality between the Lipschitz upper semicontinuity modulus and calmness moduli in linear semi-infinite inequality systems and fractured, non-convex sub-level sets. In all these cases, the topological closedness of the graph, paired with local compactness, seamlessly unlocks the exact computation of the semilocal modulus via local quantifiers.

\bibliography{bibliography}
\bibliographystyle{acm}

\end{document}